\documentclass[letterpaper, 10 pt, conference]{ieeeconf}  %

\IEEEoverridecommandlockouts                              %

\usepackage{amsmath} %
\usepackage{amssymb}  %

\usepackage{amsthm}
\usepackage{enumerate}
\usepackage{xcolor}
\usepackage{mathtools,stmaryrd}
\usepackage{empheq}
\usepackage[hidelinks]{hyperref}

\newcommand{\ie}{\emph{i.e.}}
\newcommand{\eg}{\emph{e.g.}}
\newcommand{\smallconc}[2]{\begin{bsmallmatrix} #1 \\ #2 \end{bsmallmatrix}}

\newcommand{\ox}{{\mathring{x}}}

\newcommand{\ow}{{\mathring{w}}}

\newcommand{\dparen}[1]{(\!( #1 )\!)}

\newcommand{\R}{\mathbb{R}}

\newcommand{\calL}{\mathcal{L}}
\newcommand{\calM}{\mathcal{M}}

\newcommand{\calR}{\mathcal{R}}
\newcommand{\calS}{\mathcal{S}}

\newcommand{\calW}{\mathcal{W}}
\newcommand{\calX}{\mathcal{X}}

\newcommand{\bfI}{\mathbf{I}}

\newcommand{\ul}[1]{\underline{#1}}

\newcommand{\uly}{\ul{y}}

\newcommand{\ol}[1]{\overline{#1}}

\newcommand{\oly}{\ol{y}}

\theoremstyle{plain}

\newtheorem{theorem}{Theorem}
\newtheorem{corollary}{Corollary}

\theoremstyle{definition}
\newtheorem{definition}{Definition}
\newtheorem{example}{Example}

\theoremstyle{remark}
\newtheorem{remark}{Remark}

\newcommand{\param}[1]{\dparen{#1}}

\title{\LARGE \bf
Parametric Reachable Sets Via Controlled Dynamical Embeddings
}

\author{Akash Harapanahalli and Samuel Coogan$^{1}$%
\thanks{*This work was supported in part by the Air Force Office of Scientific Research under grant FA9550-23-1-0303 and the National Science Foundation under award \#2219755.}%
\thanks{$^{1}$Akash Harapanahalli and Samuel Coogan are with the School of Electrical and Computer Engineering, Georgia Institute of Technology, Atlanta, GA, 30318, USA. \{aharapan,sam.coogan\}.gatech.edu}%
}

\begin{document}

\maketitle
\thispagestyle{empty}
\pagestyle{empty}

\begin{abstract}
In this work, we propose a new framework for reachable set computation through continuous evolution of a set of parameters and offsets which define a \emph{parametope}, through the intersection of constraints. 
This results in a dynamical approach towards nonlinear reachability analysis: a single trajectory of an embedding system provides a parametope reachable set for the original system, and uncertainties are accounted for through continuous parameter evolution.
This is dual to most existing computational strategies, which define sets through some combination of generator vectors, and usually discretize the system dynamics.
We show how, under some regularity assumptions of the dynamics and the set considered, any desired parameter evolution can be accommodated as long as the offset dynamics are set accordingly, providing a virtual ``control input'' for reachable set computation.
In a special case of the theory, we demonstrate how closing the loop for the parameter dynamics using the adjoint of the linearization results in a desirable first-order cancellation of the original system dynamics.
Using interval arithmetic in JAX, we demonstrate the efficiency and utility of reachable parametope computation through two numerical examples.

\end{abstract}

\section{INTRODUCTION}

Verifying safe operation of complex control systems is of great importance, especially when applied in safety-critical domains like autonomous driving, aerospace systems, and human-robot interactions.
A common strategy for ensuring safe behavior is through reachable set overapproximation~\cite{XC-SS:22}, where safety is verified for a guaranteed overapproximation of the set of states the system can reach under uncertain initial conditions and inputs.

For linear systems, the reachable set computation problem is more or less solved. It was shown in~\cite{StarSets:2016} that only $n+1$ simulations are required to propagate any \emph{generalized star set}, which are defined using what are called \emph{generator vectors}, and can represent any region of the state space. 
In essence, a basis of generators $G$ define a set of points through some constrained combination, \eg, $\{\ox + G\alpha : P(\alpha)\}$, where $\ox$ is some point, $\alpha$ is a vector of coefficients defining a linear combination $G\alpha$, and $P$ is some predicate evaluating True or False.
For linear systems, the superposition principle allows one to simply evolve the center $\ox$ and the generator vectors in $G$ to obtain the reachable star set.

For nonlinear systems, many computational packages~\cite{CORA,JuliaReach,zonoLAB} have adopted similar strategies by propogating sets in generator form.
The crucial operation is an iterated application of a one-step reachability computation, \eg, using a first order conservative linearization like the following, 
\begin{align*}
    \calR_k = e^{A\tau} \calR_{k-1} \oplus \calL \oplus \tau\calW,
\end{align*}
where $\tau$ is some discretization step, $A$ is the linearization of the system, and the Minkowski sums $\oplus$ handle the inputs $\calW$ and errors $\calL$ between the linearization and the true nonlinear dynamics. 
Since the main computational blocks are linear mappings and Minkowski sums, these tools have developed different set representations~\cite{SetPropogation:2021}, which greatly improve the efficiency (and closure) of these operations.
Various generator set representations have been developed, which have different tradeoffs between reachable set accuracy and computational efficiency. To name a few, there are zonotopes~\cite{CORA}, constrained zonotopes~\cite{ConstrainedZono:2016}, polynomial zonotopes~\cite{SparsePolyZono}, hybrid zonotopes~\cite{HybridZono}, ellipsotopes~\cite{Ellipsotopes}, and taylor models~\cite{TaylorModels}.
Beyond these generator-based approaches, there are many other frameworks for reachable set computation.
For instance, level set methods use the Hamilton-Jacobi equations~\cite{LevelSet:2000,HJ-reach:2017} to represent the reachable set as a level set of a solution to a partial differential equation.

In contrast to generator-based approaches, one can consider a dual formulation where the reachable set is defined by a set of constraints on the space, rather than a set of generators.
Such constraint-based frameworks abstractly propagate parameters which dynamically update the constraints, rather than generator vectors for a set. 
For instance, in linear systems, evolving dual vectors using the adjoint dynamics provides supporting hyperplanes for the reachable set~\cite{Varaiya_ReachOptimalControl:2000}.
For nonlinear systems, reachable sets represented using polytopes defined as the intersection of halfspaces (H-polytopes) can be computed by evolving offsets along the various facets on the boundary of the polytope~\cite{HarwoodBarton_Polyhedral:2016}.
We recap these two approaches~\cite{Varaiya_ReachOptimalControl:2000,HarwoodBarton_Polyhedral:2016} briefly in Section~\ref{sec:exposition}.
When the halfspaces are axis aligned, the polytopes are interval sets~\cite{Jaulin:2001} which provide efficient bounds at the cost of accuracy~\cite{scott_barton_reach:2013}.
A related framework is mixed monotonicity~\cite{CooganArcak_MMsys:2015}, which embeds the original system into a monotone dynamical system in twice the states whose evolution provides interval bounds of the reachable set~\cite{JafarpourHarapanahalliCoogan_TAC:2024}.
Using results from contraction theory~\cite{Bullo:2023}, matrix measures have also been used to propagate norm balls~\cite{MaidensArcak_matrixmeasures:2014,FanEtal_locallyopt:2016,HarapanahalliCoogan_LDI:2024}.

\paragraph*{Contributions}

In this work, we unify and generalize the previously mentioned embedding and constraint-based frameworks using the \emph{parametope}, a novel set representation for reachable set computation.
A parametope is represented as the intersection of constraints given by level sets of a parameterized function.
In Theorem~\ref{thm:param_reach}, we show how under some regularity assumptions, the parameters $\alpha$ can arbitrarily vary (continuously), as long as the offset terms $y$ compensate accordingly.
We show how the appropriate compensation can be computed as the solution to an optimization problem, or overapproximated using efficient tools like interval analysis.
In Corollary~\ref{cor:controlled_emb}, we show how to build a controlled embedding system, where a single trajectory provides overapproximating parametope sets for the original system.

Our embedding system can be viewed as a dual approach to generator representation strategies.
Rather than propagating generators, which represent real vectors in the linear case, we propagate parameters, which represent dual vectors in the linear case.
Further, rather than discretizing and handling uncertainty through Minkowski sums, we handle uncertainty dynamically, avoiding the need to compute any additional set operations.
In Section~\ref{sec:adjoint_emb}, we work out two special cases of the theory using the adjoint dynamics of the linearized system to update the parameters, and provide two numerical experiments demonstrating the computational efficiency when using interval methods in JAX.

\subsection{Notations}

Let $\leq$ denote the elementwise order on $\R^n$ satisfying $x\leq y$ if and only if $x_i\leq y_i$ for every $i=1,\dots,n$.
Let $\|\cdot\|$ denote an arbitrary norm unless otherwise specified.
Let $\|\cdot\|_{\text{op}}$ denote the operator or induced matrix norm.
For a mapping $f:\R^n\to\R^m$, let $Df:\R^n\to\R^{m\times n}$ be the Jacobian matrix $Df(x) = \frac{\partial f}{\partial x}(x)$.
Let $\bfI_n$ denote the $n\times n$ identity matrix.

\section{EXPOSITION: H-POLYTOPE EVOLUTION} \label{sec:exposition}

In this section, we briefly recall two existing results for reachable set computation using H-rep polytopes. 
For simplicity, we omit the disturbance input for these examples and focus on intuition and clarity rather than full rigor.

\subsection{Flowing Dual Vectors for Linear Systems} \label{exposition:adjoint}

Consider the linear autonomous system
\begin{align*}
    \dot{x} = Ax, \quad x(0) = x_0\in\{x : \alpha_0 (x - \ox_0) \leq y_0\},
\end{align*}
with an initial condition inside a halfspace determined by dual vector $\alpha_0$, centering point $\ox_0$, and offset $y_0$. We want to find a curve $t\mapsto\alpha(t)$, such that the trajectory $x(t) \in \{x : \alpha(t)^T (x - \ox(t)) \leq y_0\}$, for fixed offset $y_0$, and $t\mapsto \ox(t)$ as the trajectory from $\ox_0$.
As an \emph{ansatz}, suppose the desired curve $\alpha(t)$ is differentiable, then
\begin{align*}
    \frac{d}{dt} (\alpha^T (x - \ox)) 
    &= \dot{\alpha}^T (x - \ox) + \alpha^T A (x - \ox).
\end{align*}
The choice $\dot{\alpha} = -A^T \alpha$ sets the RHS to $0$ uniformly for any $x$.
These dynamics on $\alpha$ are called the adjoint dynamics, which specifically preserve the dual pairing $\alpha^Tx$ along any trajectory of the system.
Thus, if we build the following dynamics on $\alpha$ and $\ox$,
\begin{align*}
    \dot{\ox} = A\ox, \quad \dot{\alpha} = -A^T \alpha, 
\end{align*}
we have that for every $t \geq 0$,
\begin{align*}
    x(t) \in \{x : \alpha(t)^T (x - \ox(t)) \leq y_0\},
\end{align*}
since $\alpha(t)^T (x(t) - \ox(t)) \equiv \alpha_0^T (x_0 - \ox_0) \leq y_0$ is a constant function using the adjoint dynamics for $\alpha$.
For an H-rep polytope represented as $K$ halfspaces, $\{\alpha_0 (x - \ox(t)) \leq y_0\}$, for $\alpha\in\R^{K\times n}$, $y\in\R^K$, we can simultaneously flow each row $\alpha_k^T$ using the adjoint equation as $\dot{\alpha} = -\alpha A$ to obtain $x(t)\in\{\alpha(t)(x - \ox(t)) \leq y_0\}$.
This idea is the basic premise explored in~\cite{Varaiya_ReachOptimalControl:2000}.

\subsection{Fixed Dual Vectors for Nonlinear System} \label{exposition:offsets}

Another approach for flowing H-rep polytopes is to fix the halfspaces and bound the vector field along each face of the polytope. Consider the system
\begin{align*}
    \dot{x} = f(x), \quad x_0\in\{x : \alpha (x - \ox_0) \leq y_0\},
\end{align*}
with initial condition inside a compact polytope. It is shown in~\cite[Theorem 1]{HarwoodBarton_Polyhedral:2016} that under some regularity conditions, if some absolutely continuous curve $t\mapsto y(t)$ satisfies 
\begin{align*}
    \dot{y}_k(t) \geq \alpha_k^T (f(x) - f(\ox))
\end{align*}
for every $x$ along the $k$-th face of the polytope $\{\alpha(x - \ox) \leq y(t),\,\alpha_k^T(x - \ox) = y_k(t)\}$, then trajectories are in the H-polytope determined using the dual vectors $\alpha$ and the offset $y(t)$.
In other words, if we build the following dynamics on $\ox,\alpha,y$,
\begin{align*}
    \dot{\ox} = f(\ox), \quad \dot{\alpha} = 0, \quad \dot{y}_k \geq \sup_{\substack{x : \alpha(x - \ox) \leq y \\ \alpha_k^T(x - \ox) = y_k}}\alpha_k^T (f(x) - f(\ox)),
\end{align*}
we have that for every $t\geq 0$,
\begin{align*}
    x(t) \in \{x : \alpha (x - \ox(t)) \leq y(t)\},
\end{align*}
with fixed matrix $\alpha$.

As a final remark, a special case of the fixed dual vector strategy arises when the parameters are structured as $\alpha = [-\bfI_n\ \bfI_n]^T$, which given the offset $y = [-\uly^T\ \oly^T]^T$, results in the interval set $\{x\in\R^n : \uly \leq x \leq \oly\}$. 
This corresponds to setting the offset dynamics to the following,
\begin{align*}
    -\dot{\uly}_i &\geq - (\sup -e_i^T f(x)) \implies \dot{\uly}_i \leq \inf f_i (x), \\
    \dot{\oly}_i &\geq \sup e_i^T f(x) = \sup f_i (x),
\end{align*}
where the $\inf$ and $\sup$ are taken over the corresponding lower and upper $i$-th faces of $[\uly,\oly]$ ($x\in[\uly,\oly]$, $x_i=\uly_i$ or $x_i=\oly_i$),
which are easily under and over approximated using interval analysis as~\cite{scott_barton_reach:2013,JafarpourHarapanahalliCoogan_TAC:2024}.
As noted in \cite{JafarpourHarapanahalliCoogan_TAC:2024}, interval reachability connects directly to mixed monotone embedding systems \cite{CooganArcak_MMsys:2015}, which are a special case of the embedding system we construct in Corollary~\ref{cor:controlled_emb} in the next section ($U = 0$).

\section{GENERAL THEORY: PARAMETRIC REACHABLE SETS}

In the previous section, we recalled two approaches for tracking the evolution of polytopic reachable sets in H-rep.
Motivated by these examples, in this section, we pose the following questions: (i) can we bridge these theories in the nonlinear case, allowing us to dynamically update parameters constraining a set like the adjoint dynamics in the linear case, while properly handling their error using offset dynamics as in the nonlinear case? (ii) can we generalize beyond H-rep polytopes to improve the generality of the result?

\subsection{Parametopes: A New Set Representation}

We first introduce a new set representation called a \emph{parametope}, which is defined as the intersection of several parameterized sublevel sets with specified offset.

\begin{definition}
A \emph{parametope} is the set
\begin{align*}
    \param{\ox,\alpha,y}_{g} := \{x\in\R^n : g(\alpha,x - \ox) \leq y\},
\end{align*}
where $\ox\in\R^n$ (\emph{center}), $\alpha\in\R^p$ (\emph{parameters}), $y\in\R^K$ (\emph{offset}), $g:\R^m\times\R^n\to\R^K$ (\emph{nonlinearity}).
Notationally, set the \emph{$k$-th facet} of the parametope
\begin{align*}
    &\param{\ox,\alpha,y}^k_g := \param{\ox,\alpha,y}_g \cap \{x\in\R^n : g_k(\alpha,x - \ox) = y_k\}.
\end{align*}
\end{definition}

\begin{example} \label{ex:parametopes}
The following are examples of different types of sets which can be represented using parametopes.
\begin{enumerate}[i.]
    \item \label{ex:pt:H-rep} (H-rep Polytopes) Let $\alpha\in\R^{K\times n}$ and $y\in\R^K$. Then the H-rep polytope $\{x : \alpha(x - \ox) \leq y\}$ is clearly a parametope taking $g(\alpha,x - \ox) = \alpha(x - \ox)$. 
    Analyzing further, we see that for each $k$, $g_k(\alpha,x - \ox) = \alpha_k^T (x - \ox)$ is the canonical pairing between dual vector (row) $\alpha_k$ and $x - \ox$.
    \item \label{ex:pt:ellipsoid} (Ellipsoids) Let $g(\alpha,x) = x^T\alpha^T\alpha x$, $\alpha$ be a square invertible matrix, and $y > 0$. The parametope $\param{\ox,\alpha,y}_{g}$ is the ellipsoid around the center $\ox$,
    \begin{align*}
        \{x\in\R^n : (x - \ox)^T \alpha^T\alpha (x - \ox) \leq y\}.
    \end{align*}
    \item \label{ex:pt:annulus} (Annulus)
    Let $\alpha$ be a square invertible matrix, and $\|\cdot\|$ be a norm on $\R^n$, $g(\alpha,x) = \|\alpha x\|$, $0 < r_l < r_u$. Then $\param{\ox,\alpha,(-r_l,r_u)}_{(-g,g)}$ represents the annulus
    \begin{align*}
        \{x\in\R^n : r_l \leq \|\alpha(x - \ox)\| \leq r_u\}.
    \end{align*}
\end{enumerate}
\end{example}

\begin{remark}
The following remarks comment on the parametope set representation.
\begin{itemize}
\item The center $\ox$ is not necessarily inside the parametope. For instance, in Example~\ref{ex:parametopes}~\eqref{ex:pt:annulus}, $\|\alpha(\ox - \ox)\| = 0 < r_l$ so $\ox\notin \param{\ox,\alpha,(r_l,r_u)}_{(-g,g)}$.
\item A parametope is not necessarily convex or star-convex. The annulus also demonstrates this. However, structural properties of the nonlinearity and the offsets can result in additional structure of the parametope. For instance, if $g$ is a quasiconvex function in $x$ (\ie, $g(\alpha, s x + (1 - s)y) \leq \max\{g(\alpha,x),g(\alpha,y)\}$, $\forall s\in[0,1]$), the parametope $\param{\ox,\alpha,y}_g$ is convex, as the sublevel set of the quasiconvex function $x\mapsto g(\alpha,x - \ox)$.
\item If $K=1$ and $x\mapsto g(\alpha,x - \ox)$ is a positive definite, radially unbounded function, then $\param{\ox,\alpha,y}_g = \{x : g(\alpha,x - \ox) \leq y\}$ is a compact set as the sublevel set of a Lyapunov candidate.
\item Under no regularity assumptions for $g$, the parametope can represent any set $\calS$ using its indicator function, \eg, taking $g(\alpha,x) = -\mathbf{1}_\calS(x)$ and $y = -1$. However, in Theorem~\ref{thm:param_reach}, we will assume $g$ is $C^1$ with Lipschitz partial derivatives.
\end{itemize}
\end{remark}

\subsection{Reachable Parametopes}

We now present the main Theorem of this paper.
Consider a general nonlinear system defined by 
\begin{align} \label{eq:nlsys}
    \dot{x}(t) = f(x(t),w(t)), \quad x(t_0) = x_0,
\end{align}
where $x\in\R^n$ is the state of the system, $w\in\R^m$ is a disturbance input, and $f:\R^n\times\R^m\to\R^n$ is locally Lipschitz in $x$ locally uniformly on $w$. 
This means that for any compact sets $\calX\subseteq\R^n$ and $\calW\subseteq\R^m$, there exists $L_f>0$ such that $\|f(x_1,w) - f(x_2,w)\| \leq L_f\|x_1 - x_2\|$ for every $x_1,x_2\in \calX$, $w\in \calW$.
Under these assumptions, the system has a unique trajectory for some neighborhood of $t_0$ under essentially bounded disturbance $t\mapsto w(t)$, which we will denote $t\mapsto\phi_f(t,t_0,x_0,w)$.

For the proof of Theorem~\ref{thm:param_reach}, we require some additional regularity of the sets considered.

\begin{definition}[Locally Lipschitz facet {\cite[Lemma 6.1]{Still_ParametricOptLectures:2018}}] \label{def:lip_kfacets}
We say that the $k$-th facet $\param{\ox,\alpha,y}_g^k$ is locally Lipschitz in $y$ if for every $x\in \param{\ox,\alpha,y}_g^k$, there exist neighborhoods $B_{\delta}(x)$, $B_\varepsilon(y)$ and $L_k>0$ such that for every $y_1,y_2\in B_\varepsilon(y)$, the following holds:
for any $x_1\in B_\delta(x)\cap\param{\ox,\alpha,y_1}_g^k$ there exists $x_2\in B_{2\delta}(x)\cap\param{\ox,\alpha,y_2}_g^k$ such that
\begin{align*}
    \|x_1 - x_2\| \leq L_k \|y_1 - y_2\|.
\end{align*}
\end{definition}

\begin{remark}[Constraint Qualifications]
For sets defined structurally as a parameterized level set, there are some sufficient conditions, known as constraint qualifications, which ensure Lipschitz behavior with respect to the parameters.
For instance, the Linear Independence Constraint Qualification (LICQ) implies that if for every $x\in\param{\ox,\alpha,y}_g$, the vectors
\begin{align*}
    \left\{\frac{\partial g_k}{\partial x}(\alpha,x - \ox) : \forall k \text{ s.t. } g_k(\alpha, x- \ox) = y_k\right\}
\end{align*}
are linearly independent, then each facet $\param{\ox,\alpha,y}_g^k$ is locally Lipschitz as Definition~\ref{def:lip_kfacets}~\cite[Lemma 6.1]{Still_ParametricOptLectures:2018}.
Applied to the parametope, LICQ can be viewed geometrically as requiring each facet of the parametope to intersect transversally.
\end{remark}

\begin{theorem}[Parametric reachable sets] \label{thm:param_reach}
Consider the system $\dot{x} = f(x,w)$ from~\eqref{eq:nlsys}, let $g:\R^p\times\R^n\to\R^K$ be a $C^1$ map with locally Lipschitz partial derivatives, $\calW\subseteq\R^m$ be a compact set, $t\mapsto\ow(t)\in\calW$ measurable, and $t\mapsto\ox(t) = \phi_f(t,t_0,\ox_0,\ow)$.
If $t\mapsto\alpha(t),y(t)$ are absolutely continuous curves satisfying:
\begin{enumerate}[i)]
    \item for a.e. $t\in[t_0,t_f]$, $\|\dot{\alpha}(t)\| \leq A < \infty$ for some $A>0$;
    \item \label{hyp:local_lip} for every $t\in[t_0,t_f]$ and every $k\in\{1,\dots,K\}$, $\param{\ox(t),\alpha(t),y(t)}_g^k$ is locally Lipschitz as Definition~\ref{def:lip_kfacets};
    \item \label{hyp:xi_k} for a.e. $t\in[t_0,t_f]$ and every $k\in\{1,\dots,K\}$, every $x\in\param{\ox(t),\alpha(t),y(t)}_g^k$ and $w\in\calW$ satisfies $\xi_k(t,x,w) \leq \dot{y}_k(t)$, where
    \begin{align} \label{eq:xi_k}
    \begin{aligned}
        &\xi_k(t,x,w) := \frac{\partial g_k}{\partial\alpha}(\alpha(t),x - \ox(t)) [\dot{\alpha}(t)] \\
        &+ \frac{\partial g_k}{\partial x}(\alpha(t), x - \ox(t)) [f(x,w) - f(\ox(t),\ow(t))]; 
    \end{aligned}
    \end{align}
\end{enumerate}
then for any $x_0\in\param{\ox(t_0),\alpha(t_0),y(t_0)}_g$ and measurable $t\mapsto w(t) \in \calW$,  for every $t\in[t_0,t_f]$,
\begin{align*}
    \phi_f(t,t_0,x_0,w) \in \param{\ox(t),\alpha(t),y(t)}_g.
\end{align*}
\end{theorem}
The term $\xi_k$ from~\eqref{eq:xi_k} characterizes the infinitesimal change in offset for any $x$ along the $k$-th facet, due to changes in the parameters (first term) and the flow of the dynamics in the direction of the normal vector to the $k$-th facet at $x$ (second term).

\begin{proof}
Set the notation $x(t) := \phi_f(t,t_0,x_0,w)$ and $h(t) := g(\alpha(t),x(t) - \ox(t))$.
Let $T := \{t\in [t_0,t_f] : \exists i,\, h_i(t) > y_i(t)\}$, and for contradiction, suppose $T\neq\emptyset$.

Set $t_1 := \inf T$.
Using hypothesis~\eqref{hyp:local_lip}, let $\delta,\varepsilon>0$ such that $B_\delta(x(t_1))$ and $B_\varepsilon(y(t_1))$ are neighborhoods satisfying Definition~\ref{def:lip_kfacets} for every $k$, and let $L=\max\{L_1,\dots,L_K\}>0$.
By continuity, $\exists t_4>t_1$ such that $\|x(t) - x(t_1)\| < \delta/2$, $\|y(t) - y(t_1)\|_\infty < \varepsilon/2$ for every $t\in[t_1,t_4]$. 

Let $L_f>0$ be a Lipschitz constant of $f$ and $F>0$ bound $\|f(x,w)\|$ on $\ol{B_\delta(x(t_0))}\times\calW$.
Let $L_g,L_{\partial_\alpha g},L_{\partial_x g} > 0$ be Lipschitz constants for $g$, $\frac{\partial g}{\partial \alpha}$, and $\frac{\partial g}{\partial x}$ respectively, on $\{\alpha(t) : t\in[t_1,t_4]\}\times\{x - \ox(t) : x\in \ol{B_{2\delta}(x(t_1))},t\in[t_1,t_4]\}$. 

Since $t_1$ is the $\inf$ of $T$, there is a $\gamma>0$ and a $t\in(t_1,t_4]$ such that $h_i(t) > y_i(t) + \gamma$ for some $i$.
\cite[Lemma 4]{scott_barton_reach:2013} implies the existence of absolutely continuous, nondecreasing $\rho:[t_1,t_4]\to\R$ satisfying (a) $0 < \rho(t) \leq \min\{\gamma,\varepsilon/2, \delta/(2L)\}$ for every $t\in[t_1,t_4]$; and (b) $\dot{\rho}(t) \geq C \rho(t)$ for a.e. $t\in[t_1,t_4]$, where $C := L(L_{\partial_\alpha g} A + L_gL_f + 2L_{\partial_x g}F)$.
Let $t_3 := \inf\{s\in[t_1,t_4] : \exists i, h_i(s) \geq y_i(s) + \rho(s)\}$ (there exists a point satisfying this property since $\rho\leq\gamma$).
Since $\rho(t_1) > 0$, $h(t_1) < y(t_1) + \rho(t_1) \mathbf{1}$. 
Continuity implies $t_3>t_1$ and $h(t) < y(t) + \rho(t)\mathbf{1}$ for all $t\in[t_1,t_3)$.
Also, $\exists k$ s.t. $h_k(t_3) = y_k(t_3) + \rho(t_3)$, by continuity and $t_3$ being the $\inf$.
Finally, set $t_2 := \sup\{s\in[t_1,t_3] : h_k(s) \leq y_k(s)\}$. Continuity and $\rho>0$ imply $t_2 < t_3$ and $h_k(t_2) = y_k(t_2)$.

Set $t\mapsto\Tilde{y}(t)$ such that $\Tilde{y}_i(t) = \max\{h_i(t),y_i(t)\}$. Then for every $t\in[t_2,t_3]$, $h(t) \leq \Tilde{y}(t)$ and $h_k(t) = \Tilde{y}_k(t)$, thus, $x(t)\in\param{\ox(t),\alpha(t),\Tilde{y}(t)}_{g}^k$.
Since $h(t) \leq y(t) + \rho(t)\mathbf{1}$, $\|\Tilde{y}(t) - y(t)\|_\infty \leq \rho(t) \leq \varepsilon/2$ for every $t\in[t_2,t_3]$. 
Thus
\begin{align*}
    \|\Tilde{y}(t) - y(t_1)\|_\infty &\leq \|\Tilde{y}(t) - y(t)\|_\infty + \|y(t) - y(t_1)\|_\infty \\
    &\leq \varepsilon/2 + \varepsilon/2 = \varepsilon,
\end{align*}
and $y(t),\Tilde{y}(t)\in B_\varepsilon(y(t_1))$ for every $t\in[t_2,t_3]$. 
Since $x(t)\in B_\varepsilon(x(t_1))\cap\param{\ox(t),\alpha(t),\Tilde{y}(t)}_{g}^k$,
hypothesis~\eqref{hyp:local_lip} implies the existence of $z(t)\in \param{\ox(t),\alpha(t),y(t)}_g^k$ such that
\begin{align*}
    \|x(t) - z(t)\| \leq L \|\Tilde{y}(t) - y(t)\|_\infty \leq L \rho(t) \leq \delta/2,
\end{align*}
thus, triangle inequality implies
\begin{align*}
    \|z(t) - x(t_1)\| &\leq \|z(t) - x(t)\| + \|x(t) - x(t_1)\| \\
    &\leq \delta/2 + \delta/2 = \delta,
\end{align*}
so $z(t)\in B_\delta(x(t_1))$ for every $t\in[t_2,t_3]$. 

For $\zeta\in\{x,z\}$, fix the notation $\partial_\alpha g_k(\zeta) := \frac{\partial g_k}{\partial \alpha}(\alpha(t),\zeta(t) - \ox(t))$ and $\partial_x g_k(\zeta) := \frac{\partial g_k}{\partial x}(\alpha(t),\zeta(t) - \ox(t))$.
For a.e. $t\in[t_2,t_3]$, using the chain rule
\begin{align*}
    &\dot{h}_k(t)  \\
    &= \partial_\alpha g_k (x) [\dot{\alpha}(t)] + \partial_x g_k(x) [f(x(t),w(t)) - f(\ox(t),\ow(t))] \\
    &= \partial_\alpha g_k (z) [\dot{\alpha}(t)] + \partial_x g_k(z) [f(z(t),w(t)) - f(\ox(t),\ow(t))] \\
    &\quad + (\partial_\alpha g_k (x) - \partial_\alpha g_k(z)) [\dot{\alpha}(t)] \\
    &\quad + \partial_x g_k(x) [f(x(t),w(t)) - f(z(t),w(t))]  \\
    &\quad + (\partial_x g_k(x) - \partial_x g_k(z)) [f(z(t),w(t)) - f(\ox(t),\ow(t))] \\
    &\leq \xi_k(t,z(t),w(t)) \\
    &\quad + (L_{\partial_\alpha g} A + L_gL_f + 2L_{\partial_x g}F) \|x(t) - z(t)\|  \\
    &\leq \dot{y}_k(t) + (L_{\partial_\alpha g} A + L_gL_f + 2L_{\partial_x g}F) L\rho(t) \\
    &= \dot{y}_k(t) + C \rho(t) ,
\end{align*}
using hypothesis~\eqref{hyp:xi_k} since $z(t)\in\param{\ox(t),\alpha(t),y(t)}_g^k$ and $w(t)\in\calW$.
Since $\dot{\rho}(t) \geq C\rho(t)$, for a.e. $t\in[t_2,t_3]$,
\begin{align*}
    \dot{h}_k(t) \leq \dot{y}_k(t) + \dot{\rho}(t) \implies \tfrac{d}{dt} [h_k(t) - y_k(t) - \rho(t)] \leq 0,
\end{align*}
so $t\mapsto h_k(t) - y_k(t) - \rho(t)$ is nonincreasing on $[t_2,t_3]$, and
\begin{align*}
    h_k(t_2) - y_k(t_2) - \rho(t_2) \geq h_k(t_3) - y_k(t_3) - \rho(t_3).
\end{align*}
Since $h_k(t_2) = y_k(t_2)$ and $h_k(t_3) = y_k(t_3) + \rho(t_3)$, this implies $\rho(t_2) \leq 0$, a contradiction, so $T = \emptyset$.
\end{proof}

Our proof structure resembles that of~\cite[Theorem 1]{HarwoodBarton_Polyhedral:2016}, which proves a similar result for H-polytopes with fixed halfspaces.
Theorem~\ref{thm:param_reach} generalizes that result to any parametope, with time-varying parameters $\alpha$.

\begin{remark}[Connections to the literature] \label{remark:connections}
Theorem~\ref{thm:param_reach} recovers several results from the literature as special cases.
\begin{enumerate}[i)]
    \item \label{remark:connections:lin_adjoint} (Adjoint evolution) For the linear system $\dot{x} = Ax$, and polytope sets corresponding to functionals $g_k(\alpha,x - \ox) = \alpha_k^T (x - \ox)$, we get $\xi_k(t,x) = (x - \ox)^T \dot{\alpha}_k + \alpha_k^T (\dot{x} - \dot{\ox}) = - (x - \ox)^TA^T\alpha_k + \alpha_k^T A(x - \ox) = 0$. Thus, $\dot{y} = 0$ (or $y(t) \equiv y(0)$) satisfies the hypotheses of Theorem~\ref{thm:param_reach}, and recovers the result from Section~\ref{exposition:adjoint}.
    \item (Fixed dual vectors) Similarly, H-polytopes for the nonlinear system $\dot{x} = f(x,w)$ with $\dot{\alpha} = 0$ results in $\xi_k(t,x,w) = \alpha^T_k (f(x,w) - f(\ox,w))$, recovering the result from Section~\ref{exposition:offsets}.
    \item (Contraction theory) Consider the autonomous system $\dot{x} = f(x)$ for $C^1$ differentiable $f$, $\alpha$ a square invertible matrix, $\dot{\alpha} = 0$, and $g(\alpha,x - \ox) = \|x - \ox\|_\alpha = \|\alpha(x - \ox)\|$ for a differentiable norm. Using the curve norm derivative formula~\cite[Def. 17]{DavydovJafarpourBullo_NonEuc:2022}, for any $x\in\param{\ox,\alpha,y}_g^1 = \{x : \|x - \ox\|_\alpha = y\}$,
    using a weak pairing $\operatorname{WP}_\alpha$ compatible with $\|\cdot\|_\alpha$~\cite[Def. 15]{DavydovJafarpourBullo_NonEuc:2022},
    \begin{align*}
        \dot{y} &\leq \xi(t,x,w) = \frac{d}{dt} \|x - \ox\|_{\alpha} \\
        &= \frac{\operatorname{WP}_{\alpha}(f(x) - f(\ox(t)),x-\ox(t))}{\|x - \ox(t)\|_{\alpha}} \\
        &= \|x - \ox(t)\|_\alpha \sup_{s\in[0,1]} \mu(Df(sx + (1-s)\ox(t))) \leq cy
    \end{align*}
    using~\cite[Rem. 28]{DavydovJafarpourBullo_NonEuc:2022},
    where $\mu(A) = \lim_{h\searrow0} \frac{\|I + hA\|_{\text{op}} - 1}{h}$ is the logarithmic norm, and $c$ is an upper bound for $\mu(Df(x))$ on the set of interest.
    Gr\"onwall's lemma results in the bound $\|x(t) - \ox(t)\| = y(t) \leq e^{ct} y(t) = e^{ct} \|x(0) - \ox(0)\|$. 
    This recovers the results from~\cite{MaidensArcak_matrixmeasures:2014,FanEtal_locallyopt:2016,HarapanahalliCoogan_LDI:2024} which use matrix measures for norm ball reachable set computation.
\end{enumerate}
\end{remark}

\begin{remark}
For simplicity, we assumed $\calW$ was a constant compact disturbance set. The same theory allows for a time-varying disturbance set as long as $\calW(t)$ is compact.
\end{remark}

\subsection{Controlled Embeddings}

Theorem~\ref{thm:param_reach} provides a very general result for computing reachable parametopes for nonlinear systems.
In practice, the following Corollary demonstrates how the curves $t\mapsto\ox(t),\alpha(t),y(t)$ can be obtained as the solution to a controlled \emph{embedding system}.

\begin{corollary}[Controlled embeddings] \label{cor:controlled_emb}
Consider the setting of Theorem~\ref{thm:param_reach}.
Define the embedding system
\begin{align*}
    \dot{\ox} &= f(\ox,\ow), \quad
    \dot{\alpha} = U(\ox,\alpha,y),  \quad
    \dot{y} = E(\ox,\alpha,y),
\end{align*}
evolving on $\R^n\times\R^p\times\R^K$,
where $U$ is a locally Lipschitz map, and $E$ is a locally Lipschitz map satisfying
\begin{gather*}
(x,w)\in\param{\ox,\alpha,y}_g^k\times\calW \implies \\ %
\begin{aligned}
    E_k(\ox,\alpha,y) \geq &\ \frac{\partial g_k}{\partial\alpha}(\alpha,x - \ox) [U(\ox,\alpha,y)] \\
        &+ \frac{\partial g_k}{\partial x}(\alpha, x - \ox) [f(x,w) - f(\ox,\ow)],
\end{aligned}
\end{gather*}
for every $k\in\{1,\dots,K\}$.
If for every $t\in[t_0,t_f]$, $\param{\phi_{(f,U,E)}(t,t_0,(\ox_0,\alpha_0,y_0),\ow)}_g$ has locally Lipschitz $k$-facets as Definition~\ref{def:lip_kfacets},
then for any initial condition $x_0\in\param{\ox_0,\alpha_0,y_0}_g$ and disturbance map $t\mapsto w(t)\in\calW$,
\begin{align*}
    \phi_f(t,t_0,x_0,w) \in \param{\phi_{(f,U,E)}(t,t_0,(\ox_0,\alpha_0,y_0),\ow)}_g,
\end{align*}
for every $t\in[t_0,t_f]$.
\end{corollary}

Corollary~\ref{cor:controlled_emb} provides a dynamical framework for reachable set computation by embedding the original system into a new system evolving on $\R^n\times\R^p\times\R^K$. 
Computationally, we only need to simulate one trajectory of this embedding system to obtain a valid reachable parametope for the original system.
Additionally, since the map $U$ is arbitrary, this provides a new virtual input to control the parameters defining the reachable set.

\begin{remark}[Monotone embedding]
Let $U(\ox,\alpha,y) = U(\ox,\alpha)$ (no dependence on the offset $y$).
Consider the \emph{tight embedding system}, constructed by taking equality for the offset dynamics in Corollary~\ref{cor:controlled_emb}, \ie, when $\param{\ox,\alpha,y}_g^k \neq \emptyset$,
\begin{align*}
    E_k(\ox,\alpha,y) =  \displaystyle\sup_{\substack{x\in\param{\ox,\alpha,y}_g^k,\\ w\in\calW}} \xi_k(t,x,w).
\end{align*}
The Lipschitz assumption for the $k$-facets from Definition~\ref{def:lip_kfacets} also implies that the solution to the optimization $E_k$ is locally Lipschitz with respect to $\alpha,y$~\cite[Lemma 6.2]{Still_ParametricOptLectures:2018}.

The tight embedding system is a monotone system~\cite{AngeliSontag_Monotone:2003} with respect to the order $(\ox^{(1)},\alpha^{(1)},y^{(1)}) \preceq (\ox^{(2)},\alpha^{(2)},y^{(2)})$ if and only if $\ox^{(1)} = \ox^{(2)}$, $\alpha^{(1)} = \alpha^{(2)}$, $y^{(1)} \leq y^{(2)}$.
\end{remark}

\section{ADJOINT EMBEDDINGS} \label{sec:adjoint_emb}

Recall the intuition from Section~\ref{exposition:adjoint} and Remark~\ref{remark:connections}~\eqref{remark:connections:lin_adjoint}: for autonomous linear dynamics and polytope sets, a suitable choice of dynamics for the parameters $\alpha_k$ (which define linear functionals on the space) is the adjoint dynamics $\dot{\alpha}_k = -A^T \alpha_k$, which set the offset dynamics to $0$. 
When studying nonlinear systems, we propose a similar strategy: to set the parameter dynamics equal to the adjoint dynamics of the linearized system around the center trajectory $\ox$.
In this section, we introduce the feedback term $U = -\alpha Df(\ox)$ and discuss two different methods to obtain rigorous bounds for the corresponding offset dynamics defined by $\xi$.
\footnote{All experiments were performed on a desktop computer running Kubuntu 22.04, 32GB RAM, Ryzen 5 5600X. The code is available at \url{https://github.com/gtfactslab/Harapanahalli_CDC2025.git}}

\subsection{Adjoint of the Linearization} \label{subsec:adjoint_of_lin}

We first work out a special case of the theory, noting this restriction is not strictly necessary but will simplify the resulting expressions.
For the dynamics, consider the case of linear dependence on the disturbance, 
\begin{align*}
    \dot{x} &= f(x) + Cw,
\end{align*}
where $f:\R^n\to\R^n$ is $C^2$ smooth, $C\in\R^{n\times m}$ fixed.
For the parametope, suppose that for each $k$, $g_k$ factors into a composition of a scalar valued function $h_k$ and matrix multiplication with the parameters $\alpha$ as,
\begin{align*}
    g_k(\alpha_k,x - \ox) &= h_k( \alpha(x - \ox)).
\end{align*}
In this special case, using the notation $\partial h_k := \frac{\partial h_k}{\partial z}(\alpha(x - \ox))$, setting $\dot{\alpha} = U$, $\xi_k$ simplifies to 
\begin{align*}
    \xi_k(t,x,w) = \partial h_k \left[ U[x - \ox] + \alpha [f(x) - f(\ox)] + C[w - \ow] \right].
\end{align*}
Using Taylor's theorem on $f$, expanding around $\ox$ and using Lagrange's mean value remainder form, we obtain
\begin{align*}
    f(x) &- f(\ox) = Df(\ox) [x - \ox] + \tfrac12 D^2 f(z)[(x - \ox)^{\otimes 2}],
\end{align*}
for some $z = sx + (1 - s)\ox$, $s\in(0,1)$, where $D^2f(z)$ is a multilinear map of partials, and $(x - \ox)^{\otimes 2} = (x - \ox) \otimes (x - \ox)$, where $\otimes$ is the tensor product.
Plugging this expansion into the above expression yields
\begin{align*}
    \xi_k(t,x,w) &= \partial h_k [ (U + \alpha Df(\ox))[x - \ox] + C[w - \ow] \\
    & \quad\quad + \tfrac12 D^2f(z)[(x - \ox)^{\otimes2}] ].
\end{align*}
The choice $U = -\alpha Df(\ox)$ cancels the expansion of $f$ to first order, leaving 
\begin{align*}
    \xi_k(t,x,w) &= \partial h_k [ C[w - \ow] + \tfrac12 D^2f(z)[(x - \ox)^{\otimes2}] ].
\end{align*}
From another perspective, the rows $\alpha_j^T$ of $\alpha$ evolve according to the linearized adjoint dynamics as $\alpha_j = - Df(\ox)^T \alpha_j$.

Tools like \verb|immrax|~\cite{immrax} with automatic differentiation capabilities can overapproximate this $\xi_k(t,x,w)$ term using interval arithmetic, as we demonstrate in the next Example.

\begin{example}[Van der Pol oscillator]
Consider the dynamics $\dot{x} = f(x)$ of the Van der Pol oscillator,
\begin{align}
    \dot{x}_1 = x_2, \quad \dot{x}_2 = -x_1 + \mu(1 - x_1^2)x_2,
\end{align}
with $\mu=0.25$.
We consider parameters $\alpha\in\R^{2\times 2}$, with $g(\alpha, x - \ox) = \smallconc{-\alpha(x - \ox)}{\alpha(x - \ox)}$, and $y = \smallconc{-\uly}{\oly}\in\R^4$ such that $\uly \leq \oly$. 
The corresponding parametope is the symmetric polytope
\begin{align*}
    \param{\ox,\alpha,y}_g := \{\uly \leq \alpha(x - \ox) \leq \oly \}.
\end{align*}
For example, we consider the initial set determined by $\ox = [-2\ 0]^T$, $\alpha_0 = \bfI_2$, $[\uly,\oly] = [-0.125,0.125]\times[-0.00125,0.00125]$, which yields the interval set $[-2.125,-1.875]\times[-0.00125,0.00125]$. 

We use the following embedding system,
\begin{align*}
    \dot{\ox} = f(\ox), \quad \dot{\alpha} = -\alpha Df(\ox), \quad \dot{y} = E(\ox,\alpha,y),
\end{align*}
where $E$ is obtained using \verb|immrax| with automatic differentiation in JAX to automatically compute an interval overapproximation of $\xi_k(t,x,w)$ along each facet $\param{\ox,\alpha,y}_g^k$ of consideration, and set $\dot{y}_k$ to the upper bound of this guaranteed overapproximation.
The embedding system trajectory $t\mapsto\param{\ox(t),\alpha(t),y(t)}_g$ is plotted in Figure~\ref{fig:vanderpol}.
We use the Tsit5 numerical integration scheme with $t_f = 2.2\pi$, and step size $\Delta t = \frac{t_f}{500}$.
After JIT compilation, the reachable set is computed in $0.00763$s.

\begin{figure}
    \centering
    \includegraphics[width=\linewidth]{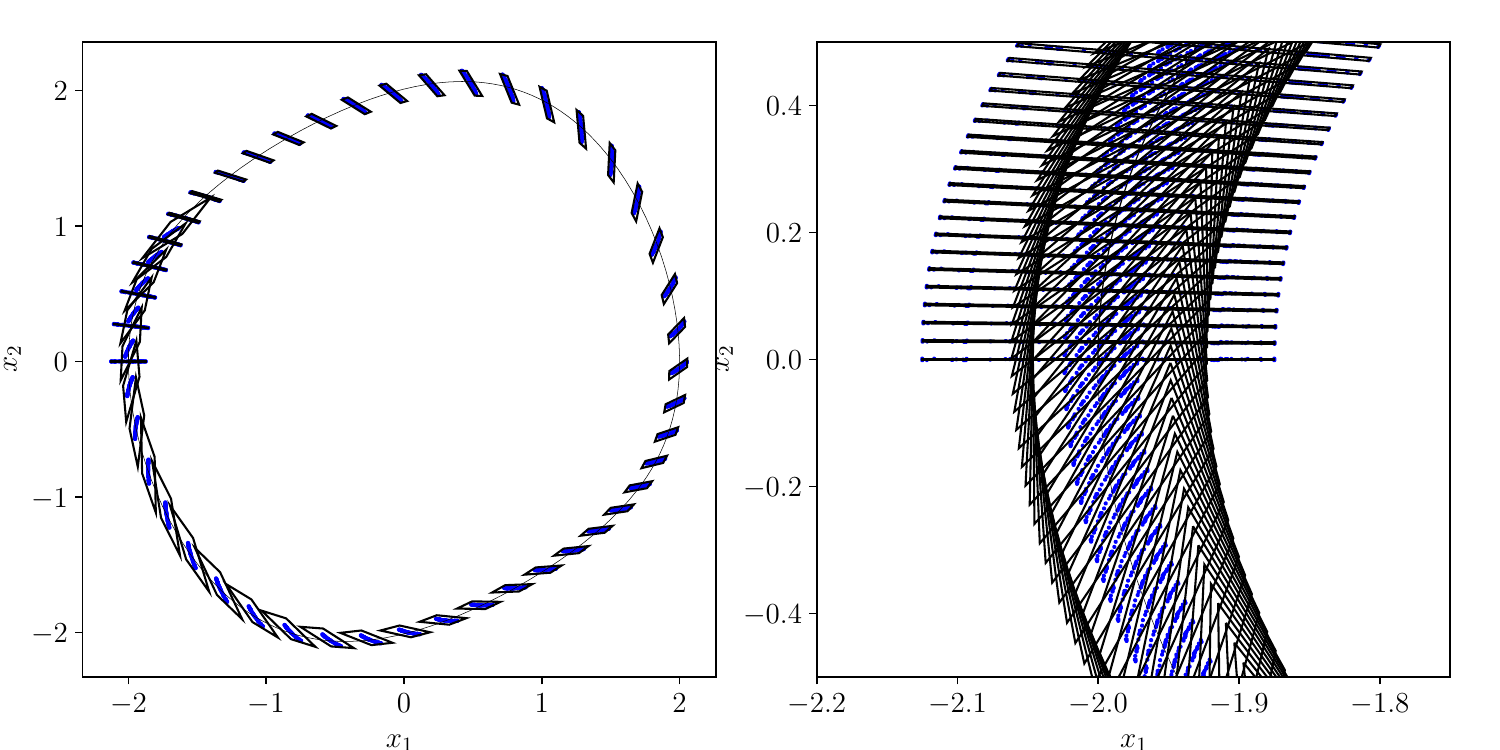}
    \caption{The overapproximating polytope reachable sets for the Van der Pol oscillator are pictured in black, with several Monte Carlo simulations in blue.
    The embedding trajectory is computed using the adjoint embedding from Section~\ref{subsec:adjoint_of_lin} from the initial set $[-2.125,-1.875]\times[-0.00125,0.00125]$. 
    \textbf{Left}: Qualitatively, the rotation of the polytope captures the spread of the Monte Carlo samples.
    \textbf{Right}: Zooming into a region around $[-2\ 0]^T$, the reachable set of the system fully passes through the initial set, verifying the stability of the limit cycle.}
    \label{fig:vanderpol}
\end{figure}
\end{example}

\subsection{Ellipsoidal Reachable Sets}

Consider the autonomous system $\dot{x} = f(x)$, for $C^1$ map $f$.
Let $\alpha_0$ be an invertible matrix, and consider the structured nonlinearity from the previous section with $h(z) = z^Tz$, resulting in $g(\alpha, x - \ox) = \|\alpha (x - \ox)\|_2^2 = (x - \ox)^T \alpha^T \alpha (x - \ox)$.
Instead of using a Taylor expansion to second order, in this section, we overapproximate the error dynamics of the system using a linear differential inclusion as~\cite{HarapanahalliCoogan_LDI:2024}, 
\begin{align*}
    f(x) - f(\ox) \in \calM (x - \ox),
\end{align*}
where $\calM\subseteq\R^{n\times n}$ is a set of matrices. 
There are different ways to obtain a set $\calM$ that satisfies this bound.
For instance, if $Df(sx + (1-s)\ox) \in \calM$ for every $s\in[0,1]$, an application of the mean value theorem implies the inclusion~\cite[Prop. 1]{HarapanahalliCoogan_LDI:2024}.
Alternatively, we can build the mixed Jacobian LDI~\cite[Thm. 1]{HarapanahalliCoogan_LDI:2024}, which is less conservative when using interval analysis to overapproximate $\calM$ since $\ox(t)$ is known~\cite[Cor. 1]{HarapanahalliCoogan_LDI:2024}.
Plugging into the expression for $\xi$ from the previous section, since $\frac{\partial h}{\partial z} (z) = 2z^T$, and setting $\dot{\alpha} = U = -\alpha Df(\ox)$ as before,
\begin{align*}
    &\xi(t,x) = 2 (x - \ox)^T \alpha^T [U[x - \ox] + \alpha[f(x) - f(\ox)]] \\
    &\leq \sup_{M\in\calM} 2 (x - \ox)^T \alpha^T \alpha (M - Df(\ox)) (x - \ox) .
\end{align*}
Recall that when optimizing over $\xi$ for the embedding dynamics, we are interested in points $x$ along the boundary, which satisfy $y = g(\alpha,x - \ox) = (x - \ox)^T \alpha^T \alpha (x - \ox)$. Rearranging, we can build a linear matrix inequality (LMI), since for any $x\in\param{\ox,\alpha,y}_g^1$,
\begin{align*}
    &\xi(t,x) 
    \leq \sup_{M\in\calM} 2 (x - \ox)^T \alpha^T \alpha (M - Df(\ox)) (x - \ox)  \\
    &\leq cy = c(x - \ox)^T \alpha^T \alpha (x - \ox) \\
    &\iff (M - Df(\ox))^T \alpha^T\alpha + \alpha^T\alpha(M - Df(\ox)) \preceq c\alpha^T\alpha
\end{align*}
for every $M\in\calM$,
and multiplying on the left by $\alpha^{-T}$ and the right by $\alpha^{-1}$,
\begin{align*}
    (\alpha(M - Df(\ox))\alpha^{-1})^T + \alpha(M - Df(\ox))\alpha^{-1} \preceq c\bfI_n.
\end{align*}
The smallest $c$ satisfying this LMI is the largest eigenvalue of the LHS. Finally, since $\lambda_{\max}$ is convex and the LHS is linear in $M$, we get that
\begin{align} \label{eq:eig_maximize}
\begin{aligned} 
    c = \max_{i} \ \lambda_{\max} ( &(\alpha(M_i - Df(\ox))\alpha^{-1})^T \\
    &+ \alpha(M_i - Df(\ox))\alpha^{-1})
\end{aligned}
\end{align}
satisfies $\xi(t,x) \leq cy$ for every $x\in\param{\ox,\alpha,y}_g^1$,
where $\{M_i\}_i$ is a set of corners satisfying $\calM\subseteq\operatorname{co}\{M_i\}_i$.

\begin{example}[Controlled ellipsoids]
\begin{figure}
    \centering
    \includegraphics[width=\linewidth]{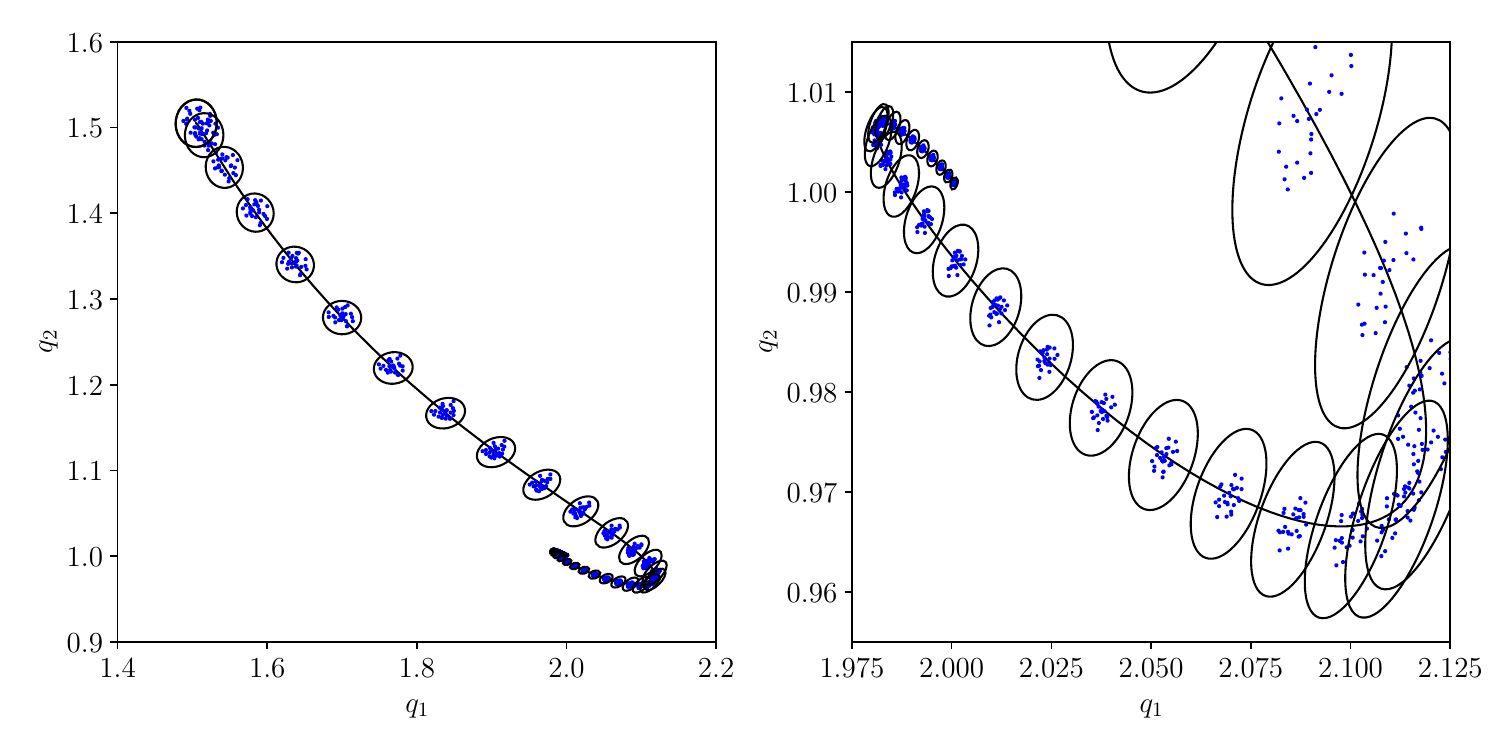}
    \caption{The projection of the overapproximating reachable ellipsoids for the robot arm onto the $q_1$-$q_2$ plane are pictured in black, around the nominal trajectory $\ox(t)$, for $t\in[0,10]$. 
    Monte Carlo samples starting in the initial ellipsoid are pictured in blue.
    \textbf{Left}: Qualitatively, the ellipsoids capture the variable spread of the Monte Carlo samples throughout the original trajectory.
    \textbf{Right}: Zooming into a region, it is clear to see that while the offset term $y$ can never decrease, the adjoint embedding dynamics results in a shrinking overapproximating reachable set.
    }
    \label{fig:robotarm}
\end{figure}
Consider the dynamics $\dot{x} = f(x)$ of a $4$ state robot arm from~\cite{AngeliSontagWang:2000,FanEtal_locallyopt:2016}
\begin{align*}
    \dot{q}_1 &= z_1, \quad \dot{q}_2 = z_2, \\
    \dot{z}_1 &= \tfrac{1}{mq_2^2 + ML^2/3}(-2mq_2z_1z_2 - k_{d_1}z_1 + k_{p_1}(u_1 - q_1)), \\
    \dot{z}_2 &= q_2 z_1^2 + \tfrac1m (- k_{d_2} z_2 + k_{p_2} (u_2 - q_2)),
\end{align*}
with $u_1 = 2$, $u_2 = 1$, $m=M=1$, $L=\sqrt{3}$, $k_{p_1} = 2$, $k_{p_2} = 1$, $k_{d_1} = 2$, $k_{d_2} = 1$.
We use the initial set $\param{\ox_0,\alpha_0,y_0}_g$, where $\alpha_0 = P^{1/2}$ with $P$ and $\ox_0$ from~\cite{FanEtal_locallyopt:2016}, and the following embedding system,
\begin{align*}
    \dot{\ox} = f(\ox), \quad \dot{\alpha} = -\alpha Df(\ox), \quad \dot{y} \leq cy,
\end{align*}
where $c$ is the solution to the eigenvalue maximization~\eqref{eq:eig_maximize} initialized as follows.
We use \verb|immrax.mjacM| to compute the interval mixed Jacobian matrix $[\calM]$ with the state variable order $[z_1\ z_2\ q_1\ q_2]^T$~\cite{immrax}, sparsely extract the $64$ corners (there are only $6$ nonconstant elements), and evaluate~\eqref{eq:eig_maximize} using \verb|jnp.linalg.inv| to find $\alpha^{-1}$ and \verb|jnp.linalg.eigvalsh| for $\lambda_{\max}$.
Comparing to the results from~\cite{FanEtal_locallyopt:2016,HarapanahalliCoogan_LDI:2024}, rather than occasionally resynthesizing a $P=\alpha^T\alpha$ matrix in a semi-definite program, which is valid for a specified horizon, we automatically and dynamically update the parameters using the adjoint equation.
We use Euler integration with $t_f = 10$, and step size $\Delta t = 0.01$.
The full embedding trajectory is computed in $3.8$s.

\end{example}

\section{CONCLUSION}

In this work, we presented a new framework for parametric reachable set computation through dynamical embeddings.
This paper, while primarily focused on the theory behind reachable set computation using parametope embeddings, demonstrated good accuracy, scalability, and computational tractability.
We believe there are many directions for future work on the computational side.
As with generator-based strategies, we suspect that enforcing additional structure to the parametope like constrained zonotopes, polynomial zonotopes, etc., will provide better capabilities for improving reachable set efficiency and accuracy.
For instance, if $g$ is defined by a polynomial whose coefficients are the parameters $\alpha$, a higher order analogue of the adjoint embedding may allow for higher order polynomial cancellations.
Additionally, the arbitrary dynamics for $\alpha$ leave room for synthesizing various control policies for the embedding system to virtually to control the growth of the reachable set.

\bibliographystyle{ieeetr}
\bibliography{citations.bib}

\end{document}